\documentstyle[amscd]{amsart}
\numberwithin{equation}{section}
\theoremstyle{plain}
   \newtheorem{thm}{Theorem}[section]
   \newtheorem{lem}[thm]{Lemma}
   \newtheorem{prop}[thm]{Proposition}
   
   \newtheorem{conj}[thm]{Conjecture}
   
   \newtheorem{claim}[thm]{Claim}
\theoremstyle{definition}
   \newtheorem{defn}[thm]{Definition}
   \newtheorem{example}[thm]{Example}
\theoremstyle{remark}
   \newtheorem{rem}[thm]{Remark}

\begin{document}
\title[Arithmetic Hodge structure]
{Arithmetic Hodge structure and higher Abel-Jacobi maps}
\author[M. Asakura]
{Masanori Asakura}
\begin{abstract}
In this paper, 
we show some applications to algebraic cycles by using
higher Abel-Jacobi
maps which were defined in [the author:
Motives and algebraic de Rham cohomology]. 
In particular, we prove that the Beilinson conjecture on algebraic 
cycles over number fields implies the Bloch conjecture on zero-cycles
on surfaces.
Moreover, we construct a zero-cycle on a product of curves
whose Mumford invariant vanishes, but not
higher Abel-Jacobi invariant.
\end{abstract}

\maketitle
%
%symbol for simplicity
\def\Spec{{\operatorname{Spec}}}
\def\Pic{{\mathrm{Pic}}}
\def\Ext{{\mathrm{Ext}}}
\def\NS{{\mathrm{NS}}}
\def\Picv{{\mathrm{Pic^0}}}
\def\Div{{\mathrm{Div}}}
\def\CH{{\mathrm{CH}}}
\def\deg{{\mathrm{deg}}}
\def\dim{{\operatorname{dim}}}
\def\codim{{\operatorname{codim}}}
\def\Coker{{\operatorname{Coker}}}
\def\ker{{\operatorname{ker}}}
\def\Image{{\operatorname{Image}}}
\def\Aut{{\mathrm{Aut}}}
\def\Hom{{\mathrm{Hom}}}
\def\Proj{{\mathrm{Proj}}}
\def\Sym{{\mathrm{Sym}}}
\def\Image{{\mathrm{Image}}}
\def\Gal{{\mathrm{Gal}}}
\def\GL{{\mathrm{GL}}}
\def\End{{\mathrm{End}}}
\def\P{{\bold P}}
\def\C{{\bold C}}
\def\Q{{\bold Q}}
\def\Z{{\bold Z}}
\def\F{{\bold F}}
\def\nb{\nabla}
\def\Lam{\Lambda}
\def\lam{\lambda}
\def\Om{\Omega}
\def\k{\kappa}
\def\l{\ell}
\def\z{\zeta}
\def\Te{\Theta}
\def\te{\theta}
\def\sg{\sigma}
\def\ve{\varepsilon}
\def\lra{\longrightarrow}
\def\ra{\rightarrow}
\def\la{\leftarrow}
\def\lla{\longleftarrow}
\def\Lra{\Longrightarrow}
\def\Lla{\Longleftarrow}
\def\da{\downarrow}
\def\hra{\hookrightarrow}
\def\lmt{\longmapsto}
\def\ot{\otimes}
\def\op{\oplus}
%decoration
\def\wt#1{\widetilde{#1}}
\def\ol#1{\overline{#1}}
\def\us#1#2{\underset{#1}{#2}}
\def\os#1#2{\overset{#1}{#2}}
\def\lim#1{\us{#1}{\varinjlim}}
%
%extra symbol
\def\et{\'{e}tale }
\def\MCb{\underline{{\mathrm{M}}}_{\C}}
\def\mot{{\cal{M}}}
\def\vg{\varGamma}
\def\HS#1{\underline{{\mathrm{HS}}}(#1)}
\def\QHS#1{\underline{{\mathrm{HS}}}(#1)_{\Q}}
\def\MHS{\underline{{\mathrm{MHS}}}}
\def\rMHS{{\mathrm{MHS}}}
\def\QMHS{\underline{{\mathrm{MHS}}}_{\Q}}
\def\MHM{{\mathrm{MHM}}}
\def\HM{{\mathrm{HM}}}
\def\VMHS{{\mathrm{VMHS}}}
\def\MM{\underline{{\mathrm{M}}}}
\def\Alb{{\mathrm{Alb}}}
\def\Perv{{\mathrm{Perv}}}
\def\gr{{\mathrm{gr}}}
\def\Gr{{\mathrm{Gr}}}
\def\cont{{\mathrm{cont}}}
\def\O{{\cal{O}}}
\def\X{{\Xi}}
\def\rank{{\mathrm{rank}}}

\def\D{{\Bbb D}}
\def\dR{{\mathrm{dR}}}
\def\crys{{\mathrm{crys}}}
\def\Gr{{\mathrm{Gr}}}
\def\dego{\deg=0}
\def\Qb{\bar{\Q}}
\def\prim{0}%{{\mathrm{prim}}}
\def\al{\alpha}
\def\fp{{\frak{p}}}
\def\fm{{\frak{m}}}
\section{Introduction}
In \cite{masanori}, we defined a certain Hodge-theoretic structure 
for an arbitrary variety $X$ over the complex number field $\C$ by using 
the theory of mixed Hodge module due to Morihiko~Saito. 
We call it an {\it arithmetic Hodge structure} of $X$. 
We defined the {\it higher Abel-Jacobi maps} 
from Bloch's higher Chow groups $\CH^r(X,m;\Q) = \CH^r(X,m)\ot\Q$ of $X$ 
to the extension groups in the category of arithmetic Hodge structures. 
(In this paper we will simply write $\CH^r(X,m)$ instead
of $\CH^r(X,m;\Q)$).

\vspace{0.1cm}

The purpose of this paper is to show some applications of
our higher Abel-Jacobi maps to algebraic cycles. 
In particular, we will give the correct proofs of 
the results announced in \cite{masanori}. 

\vspace{0.1cm}

We are interested in zero-cycles on surfaces,
in particular,
the kernel $T(X)$ of a Albanese map
$$\CH_0(X)_{\dego}\lra \Alb(X).$$
If the geometric genus $p_g$ is not zero, 
it was shown by D.~Mumford
that $T(X)$ is enormous (\cite{mumford}).
Conversely, S.~Bloch conjecture:
\begin{conj}[\cite{bloch} Lecture 1]\label{Bloch}
Let $X$ be a nonsingular projective surface with $p_g=0$.
Then $T(X)=0$.
\end{conj}
If our higher Abel-Jacobi maps for zero-cycles on surfaces
are injective, the Bloch conjecture \ref{Bloch} is true.
More generally, the injectivity of higher Abel-Jacobi maps
(Conjecture~\ref{mainconj}) 
implies the Bloch-Beilinson conjecture, that is, the existence of
motivic filtrations on Chow groups
(Proposition~\ref{mainbloch326}).

Our first main result is to reduce the Bloch conjecture \ref{Bloch}
to the Beilinson conjecture on the algebraic cycles over number fields:
\begin{thm}[Theorem~\ref{asakura2}, cf.\cite{masanori} Theorem~4.10]
\label{asakura2intro}
The Beilinson conjecture~$\ref{Beilinson}$~$(1)$, $(2)$ 
for codimension $r=2$
implies the Bloch conjecture $\ref{Bloch}$.
\end{thm} 

\vspace{0.1cm}

Another remarkable advantage of our higher Abel-Jacobi maps is that 
we can do ``explicit calculations''. 
Moreover, it gives a stronger invariant than 
higher normal functions (cf.\cite{AS}), 
although our higher Abel-Jacobi maps are defined along the idea of them. 
The reason of it is that our arithmetic Hodge structure contains a datum of 
$\Q$-coefficient perverse sheaf. Let us explain this more precisely. 
In \cite{masanori}, we defined the spaces of 
Mumford invariants $\X^{p,q}_X(r)$
and
$\Lambda^{p,q}_X(r)$, and constructed the natural maps 
$\Ext^{p}_{\MM(\C)}(\Q(0),H^{q}(X)(r))\ra
\X^{r-p,q-r+p}_X(p) \ra \Lambda^{r-p,q-r+p}_X(p)$. 
When $X$ is a surface and $p=q=r=2$, it gives the classical 
Mumford invariant (see \cite{mumford}, \cite{voisin}). 
So far, many results on algebraic cycles of codimension $r\geq2$ 
have been obtained only by observing it precisely. 
There exist, however, algebraic cycles which cannot be 
captured by only the Mumford invariants. 
Even for such a cycle, our higher Abel-Jacobi maps work quite effectively. 
The following is our second main result, which gives one example
to the above:
\begin{thm}[Theorem~\ref{asakura1}, cf.\cite{masanori} Theorem~4.5]
\label{asakura1intro}
Let $C$ be a projective nonsingular curve over $\Qb$ 
with genus $g\geq2$, such 
that $\rank \NS(C\times C)=3$.
Let $O\in C(\Qb)$ be a $\Qb$-valued point such that $K_C-(2g-2)\cdot O$
is not $\Q$-linearly equivalent to $0$.
Let $P\in C(\C)\setminus C(\Qb)$ be any $\C$-valued point which is not
$\Qb$-valued one.
Put $X:=C_{\C}\times C_{\C}$ and $z:=(P,P)-(P,O)-(O,P)+(O,O)\in T(X)$.
Then we have
$$
\xi_X^{2}(z)=0 \quad \text{in }\quad
\X^{0,2}_X(2),
$$
but
$$
\rho_X^2(z)\not=0 \quad \text{in }
\Ext_{\MM(\C)}^2(\Q(0),H^2(X)(2)).
$$
\end{thm}
For the proof of this theorem we use the modified diagonal cycle and 
through the proof of it, we see that images of algebraic cycles 
under our higher Abel-Jacobi map can be calculated 
explicitly by reducing it to the calculations of the usual Abel-Jacobi map. 

\medskip

Although M.~Green and P.~Griffiths defined 
arithmetic Hodge structure independently to the author~
(the ``arithmetic Hodge structure'' was 
named by them which seems quite nice, so the author also uses it 
in this paper), there are some differences 
in definition between theirs and ours. 
The most essential difference is that their 
definitions do not consider the datum of the $\Q$-structure which is taken 
in ours. Therefore, the category of their 
arithmetic Hodge structures does not become an abelian category, 
but only an exact category. Nevertheless, 
extension groups can be defined also in their category, 
and showed that cycle maps have non-trivial images. 
However, the conjecture \ref{mainconj} does not 
hold for their category perhaps from the lack of considering 
$\Q$-structure in definition.

\medskip

We will explain how this paper consists. 
In \S2, we review Carlson's description about extension groups
in the category of graded polarizable mixed Hodge structures,
and usual Abel-Jacobi maps.
In \S3, we review arithmetic Hodge structures, higher Abel-Jacobi
maps and Mumford invariants of algebraic cycles, which are
constructed in \cite{masanori}.
In \S4, we will prove Theorem \ref{asakura1intro}.
In \S5, we will prove Theorem \ref{asakura2intro}.

\par
\vspace{0.25cm}
\hspace{-6.0mm} {\bf Acknowledgement}\par
\vspace{0.15cm}
The author would like to express his sincere 
gratitude to Professor Shuji Saito for 
many helpful suggestions and stimulating discussions. 

He also thanks to Professor Shin-ichi 
Mochizuki for teaching him the proof of Lemma~\ref{Mochizuki}, 
and to Doctor Ken-ichiro Kimura for teaching him about Griffiths groups. 
Finally, he would like to thank to 
Professor Sampei Usui and Professor Takeshi Saito 
for fruitful discussions and unceasing encouragement. 

The author is supported by JSPS Research Fellowships for Young Scientists.
\par
\vspace{0.25cm}
\hspace{-6.0mm} {\bf Notation and Conventions}\par
\vspace{0.15cm}
\begin{enumerate}
\item
A {\it variety} means a quasi-projective algebraic variety over a field.
We mainly work with algebraically closed fields of characteristic 0
(e.g. $\C$, $\Qb$).
\item
For a variety $X$ over a field $k$, we denote $X(S)=\Hom_k(S,X)$ the set
of {\it $S$-valued points} of $X$.
\item
For a variety $X$ over $\C$, $X^{an}$ denotes the associated analytic
space: $X^{an}=X(\C)$.
\item
We denote the kernel of the cycle map $\CH^r(X)\ra H^{2r}(X)$ by
$\CH^r(X)_{\hom}$ the subgroup of homologically trivial cycles.
Here $H^{\bullet}(X)$ is a Weil cohomology (e.g. Betti cohomology,
algebraic de Rham cohomology, \et cohomology, and so on).
We also write $\CH_0(X)_{\hom}=\CH_0(X)_{\dego}$.
\item
In this paper, we fix an embedding $\Qb\hra\C$.
\end{enumerate}

%%%%%%%%%%%%%%%%%%%%%%%%%%%%%%%%%%%% Carlson %%%%%%%%%%%%%%%%%%%%%%%%%%%%%%%
\section{The Carlson isomorphism on the extensions of mixed Hodge
structures}
The category $\rMHS$ of graded polarizable $\Q$-mixed Hodge structures
is an abelian category, but not, semi-simple. 
The extension groups in $\rMHS$ have the well-known
explicit description due to J.Carlson (\cite{carlson}).

\begin{thm}
Let $H=(H_{\Q},W_{\bullet},F^{\bullet})$ be a graded polarizable 
$\Q$-mixed Hodge structure. We write $H_{\C}=H_{\Q}\ot\C$. Then
\begin{enumerate}
\renewcommand{\labelenumi}{(\theenumi)}
\item
$
\Ext_{\rMHS}^{1}(\Q(0),H)=
\Ext_{\rMHS}^{1}(\Q(0),W_0H)\simeq
W_{-1}H_{\C}/W_{-1}H_{\C}\cap(F^0W_0H_{\C}+W_0H_{\Q}).
$
Here an element $\xi \in W_{-1}H_{\C}$ corresponds to the following extension
of mixed Hodge structures:
$$
0 \lra W_0H \lra \wt{H} \lra \Q(0) \lra 0,
$$
where $\wt{H}=(\wt{H}_{\Q}, W_{\bullet}, F^{\bullet})$ is the mixed Hodge
structure with $\wt{H}_{\Q}=W_0H_{\Q}\op \Q(0)$, and the weight
filtration $W_0\wt{H}_{\Q}=\wt{H}_{\Q}$, $W_{\l}\wt{H}_{\Q}=W_{\l}H_{\Q}$
($\l \leq -1$), and the Hodge filtration 
$F^p\wt{H}_{\C}=F^pW_0H_{\C}$ ($p \geq 1$), 
$F^q\wt{H}_{\C}=F^qW_0H_{\C}+F^0\wt{H}_{\C}$ ($q \leq -1$) and 
$$
F^0\wt{H}_{\C}=F^0W_0H_{\C}+ \C \cdot (\xi , 1).
$$
\item
$\Ext_{\rMHS}^{p}(\Q(0),H)=0$ for $p \geq 2$. $($This is a formal
 consequence of the fact that the functor $\Ext^1_{\rMHS}(\Q(0), -)$ 
is right exact$)$.
\end{enumerate}
\label{carl}
\end{thm} 
\begin{rem}
\begin{enumerate}
\renewcommand{\labelenumi}{(\theenumi)}
\item
It is easy to see that $\Ext_{\rMHS}^{\bullet}(H_1,H_2)\simeq
\Ext_{\rMHS}^{\bullet}(\Q(0),H_1^*\ot H_2)$.
\item
The above description (=Theorem~\ref{carl}~(1), which I was learned
from Morihiko Saito) is different from
the one of the extension groups in the category $\MHS$ of
 mixed Hodge structures which are not necessarily graded 
polarizable. Originally, J.~Carlson calculated the extensions
in $\MHS$ (\cite{carlson}):
$$
\Ext_{\MHS}^{1}(\Q(0),H)\simeq
W_{0}H_{\C}/F^0W_0H_{\C}+W_0H_{\Q}.$$
\end{enumerate}
\end{rem}

\bigskip

Let $X$ be a projective
nonsingular variety $X$ over $\C$. By Theorem \ref{carl} (1), we have
\begin{equation}\label{carliso}
\Ext_{\rMHS}^{1}(\Q(0),H^{2r-1}(X,\Q(r)))\simeq J^r(X).
\end{equation}
Here $J^r(X)=H^{2r-1}(X^{an},\C)/F^r+H^{2r-1}(X^{an},\Q(r))$ is 
the {\it $r$-th intermediate
Jacobian} of $X$ (modulo torsion).
It is isomorphic to the Picard variety
${\mathrm Pic}^0(X)(\C)$ if $r=1$, and the Albanese variety $\Alb(X)(\C)$
if $r=\dim X$.

\medskip

Let us recall the Abel-Jacobi map
\begin{equation}\label{AJ}
\rho:
\CH^r(X)_{\hom} \lra \Ext^1_{\rMHS}(\Q(0),H^{2r-1}(X,\Q(r))).
\end{equation}
Due to the formalism of mixed Hodge modules (\cite{msaito1}, 
\cite{msaito}), there is the cycle map 
$$
\CH^r(X)\lra \Ext^{2r}_{\MHM(X)}(\Q_{X}(0),\Q_{X}(r)).
$$
It induces the following commutative
diagram:
\begin{equation}\label{algab33}
\begin{matrix}
0& &0\\
\downarrow && \downarrow\\
\CH^r(X)_{\hom}&\os{\rho}{\lra}&\Ext^1_{\rMHS}(\Q(0),H^{2r-1}(X,\Q(r)))\\
\downarrow && \downarrow\\
\CH^r(X)&\lra&\Ext^{2r}_{\MHM(X)}(\Q_{X}(0),\Q_{X}(r))\\
\downarrow && \downarrow\\
\CH^r(X)/\CH^r(X)_{\hom}&\lra& H^{2r}(X,\Q)\cap H^{r,r}\\
\downarrow && \downarrow\\
0& &0.
\end{matrix}
\end{equation}
The top horizontal arrow $\rho$ is the Abel-Jacobi map.

\medskip

On the other hand, there is the classical definition of Abel-Jacobi maps
due to Griffiths and Weil.
Let $Z\in \CH^r(X)_{\hom}$ be an algebraic cycle.
There is a topological $(2\dim X-2r+1)$-cycle $\Gamma$ whose boundary
is $Z$: $\partial \Gamma =Z$.
Then the classical Abel-Jacobi class is defined as follows: 
\begin{equation}\label{class33}
\rho(Z)=\sum_{j=1}^g\left(\int_{\Gamma}\omega_j\right)\omega_j^*\in J^r(X),
\end{equation}
where $\omega_1,\cdots,\omega_g\in F^{n-r+1}H^{2n-2r+1}(X,\C)$ is a basis,
and $\omega_1^*,\cdots,\omega_g^*\in H^{2r-1}(X,\C)/F^r$ 
denotes the Serre dual
class of those: $\langle \omega_i,\omega_j^*\rangle=\delta_{ij}$.

\medskip

The following is well-known (cf. \cite{EZ}):
\begin{prop}\label{algAJ}
The classical Abel-Jacobi map \eqref{class33} coincides with
the previous one \eqref{AJ} under the Carlson isomorphism
\eqref{carliso}.\qed 
\end{prop}

\medskip

Let $k$ be an algebraically closed subfield of $\C$.
Put $\MM_{k}=\MHM(\Spec k)$ which consists of objects
$H=(H_{\Q},H_{k},F^{\bullet},W_{\Q,\bullet},W_{k,\bullet},i)$
where
\begin{itemize}
\item
$H_{\Q}$ is a finite $\Q$-vector space,
\item
$H_{k}$ is a finite $k$-vector space,
\item
$F^{\bullet}$ is a finite decreasing filtration on $H_{k}$
(called the Hodge filtration),
\item
$W_{\Q,\bullet}$ is a finite increasing filtration on $H_{\Q}$
(called the weight filtration),
\item
$W_{k,\bullet}$ is a finite increasing filtration on $H_{k}$
(called the weight filtration),
\item
$i:H_{k}\ra H_{\C}:=H_{\Q}\ot\C$ is a $k$-linear map
(called the comparison map)
such that $i_{\C}:H_{k}\ot_k\C\ra H_{\C}$ is bijective,
\end{itemize}
which satisfy
\begin{enumerate}\renewcommand{\labelenumi}{(\theenumi)}
\item
$ W_{\Q,\bullet}$ and $W_{k,\bullet}$ is compatible under
the comparison map $i$,
\item
$(H_{\Q},W_{\Q,\bullet},i(F^{\bullet}))$ is a mixed Hodge structure.
\item
There is a polarization form on each graded component
$\Gr^W_{\l}H=(\Gr^W_{\l}H_{\Q},\Gr^W_{\l}H_{k},F^{\bullet},i)$
defined over $k$.
\end{enumerate}

We have an analogous description of the extension groups in $\MM_{k}$
(the proof is similar to the one of Theorem~\ref{carl}).

\begin{prop}\label{carl325}
\begin{enumerate}
\renewcommand{\labelenumi}{(\theenumi)}
\item
\begin{align*}
\Ext^1_{\MM_{k}}(\Q(0),H)
&=\Ext^1_{\MM_{k}}(\Q(0),W_0H)\\
&=W_{-1}H_{\C}/W_{-1}H_{\C}\cap(i(F^0W_0H_{k})+W_0H_{\Q}).
\end{align*}
\item
$\Ext^{p}_{\MM_{k}}(\Q(0),H)=0$ if $p\geq2$.
\item
$\Ext^{\bullet}_{\MM_{k}}(H_1,H_2)=
\Ext^{\bullet}_{\MM_{k}}(\Q(0),H_1^*\ot H_2)$. 
\end{enumerate}
\end{prop}

%%%%%%%%%%%%%%%%%%%%%%%%%%%%%% Arithmetic Hodge structure %%%%%%%%%%%%%%%
\section{Arithmetic Hodge structure: review of \cite{masanori}}
We review the notion of arithmetic Hodge structures, 
higher Abel-Jacobi maps and the Mumford invariants
which are introduced in \cite{masanori}.

\subsection{Arithmetic Hodge structure}
Let $X$ be a quasi-projective nonsingular variety over $\C$.
Then $X$ is defined by finitely many equations which possess
finitely many coefficients.
By considering the coefficients as parameters of a space $S$,
we can obtain a model $f:X_S\ra S$ and the Cartesian diagram:
\begin{equation}\label{model}
\begin{CD}
X_S @<<< X\\
@V{f}VV  @VVV\\
S @<{a}<<\Spec \C,
\end{CD}
\end{equation}
where $S$ is a nonsingular variety over $\Qb$, and the map $a$ factors
through the generic point
$\Spec~\Qb(S) \hookrightarrow S$.
We define the abelian category
$$\MM(X)=\lim{X_S} \MHM(X_S)$$
Here $\MHM(X_S)$ is the category of mixed Hodge modules on the variety $X_S$
over $\Qb$. In the above limit,
$X_S$ runs over all models \eqref{model}, and for a morphism
$j:X_{S'}\ra X_S$ of the models, we take the pull-back 
$j^*:\MHM(X_S)\ra\MHM(X_{S'})$.
We call it the category of {\it arithmetic Hodge modules}.
In particular, we call $\MM(\C):=\MM(\Spec\C)$
the category of {\it arithmetic Hodge structures}.

\bigskip

The {\it realization functor} $r_X:\MM(X)\ra \Perv(X)$ are constructed
as follows.
In the diagram \eqref{model}, the morphism $a:\Spec~\C\ra S$ induces
the $\C$-morphism ${\cal O}_S\ot\C\ra \C$, which defines a closed point
$s\in S^{an}=S\ot_{\Qb}\C$. Let $i_s:\Spec~\C\hra S^{an}$ be the 
corresponding inclusion. Then $a$ factors as $\Spec~\C\os{i_s}{\hra} S^{an}
\ra S$, and similarly $X\os{i_{X,s}}{\hra} X_S^{an}\ra X_S$.
We define $r_{X_S}:\MHM(X_S)\ra \Perv(X)$ the functor which
maps the perverse sheaf $K^{\bullet}_{\Q}$ of a mixed Hodge module
to ${}^pH^0i_{X,s}^{*}(K^{\bullet}_{\Q})$.
These functors are compatible for any transition $X_{S'}\ra X_S$.
Therefore,
passing to the limit, we get the realization functor $r_X=\lim{X_S}~r_{X_S}$.

\vspace{0.1cm}

Although the model \eqref{model} is not uniquely determined, any two models
$X_{S_1}$ and $X_{S_2}$
can be imbbeded into the following diagram:

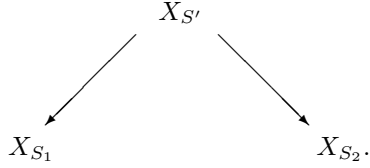
\setlength{\unitlength}{1mm}
\begin{picture}(155,35)(-76,-25)
\put(0,1){$X_{S'}$}
\put(8,-1){\vector(1,-1){12}}
\put(-3,-1){\vector(-1,-1){12}}
\put(-20,-17){$X_{S_1}$}
\put(21,-17){$X_{S_2}$.}
\label{commutativediagram}
\end{picture}
Therefore the Tate Hodge module $\Q_{X_S}(r)$ defines a well-defined
object $\Q_X(r)$ in $\MM(X)$, which we call 
the {\it arithmetic Tate Hodge module}. 

\vspace{0.1cm}

The following is straightforwards due to the formalism of
mixed Hodge modules.
\begin{prop}\label{formalism81}
\begin{enumerate}
\renewcommand{\labelenumi}{$($\theenumi$)$}
\item
There are the standard operations on the derived category
of bounded complex of arithmetic Hodge modules:
$$f_*,~f_!,~f^*,~f^!,~{\Bbb D},~\ot,~\underline{\Hom}.$$
Those functors 
satisfies the adjointness, projection formulas, and are compatible
with the ones on perverse sheaves 
under the realization functor $r_X:\MM(X)\ra \Perv(\Q_X)$.
\item
Let $f:X\ra Y$ be a proper morphism of
nonsingular varieties, and $M$ a pure object of $\MM(X)$.
Then the decomposition theorem 
of type Beilinson-Bernstein-Deligne-Gabber
$$f_*M\simeq \us{k}{\op}H^kf_*M[-k]
$$
holds.
\end{enumerate}
\end{prop}
We denote the arithmetic Hodge structure
$H^k(X,\Q(r))=H^k(X)(r):=H^kf_*\Q_X(r)$ for the structure
morphism $f:X\ra \Spec\C$.
As an immediate corollary of Proposition~\ref{formalism81}, we have
the Leray spectral sequence
\begin{equation}\label{leray81}
E_2^{pq}=\Ext^p_{\MM(Y)}(\Q_Y(0),H^qf_*\Q_X(r))\Longrightarrow
\Ext^{p+q}_{\MM(X)}(\Q_X(0),\Q_X(r)),
\end{equation}
for a morphism $f:X\ra Y$ of nonsingular varieties. 
If $f$ is proper, it degenerates 
at $E_2$-terms.

\subsection{The space of Mumford invariants}
Recall the spaces of Mumford invariants
$\X_X^{p,q}(r)$ and $\Lambda_X^{p,q}(r)$ for a projective
nonsingular variety $X$ over $\C$ (see for details, 
\cite{masanori} \S3).

Let $H^{q}_{\dR}(X/\C):=H^{q}(X,\Omega^{\bullet}_{X/\C})$
be the algebraic de Rham cohomology of $X$, and
$F^{r}H^{q}_{\dR}(X/\C):=H^{q}(X,\Omega^{\bullet\geq r}_{X/\C})$
the Hodge filtration.
There is the Gauss-Manin connection
$$
\nb :H_{\dR}^q(X/\C)\lra H_{\dR}^q(X/\C)\ot_{\C}\Omega^1_{\C/\Qb}.
$$
It satisfies the Griffiths transversality
$
\nb(F^{r})\subset 
F^{r-1}\ot\Omega^1_{\C/\Qb}$.
Then we define the space of the {\it Mumford invariants}
$\X_X^{p,q}(r)$ $($resp. $\Lambda_X^{p,q}(r)$ $)$
as the cohomology at the middle term
of the following complex:
$$
F^{p+1}H^{p+q}_{\dR}(X/\C)\ot \Omega^{r-1}_{\C/\Qb}
\os{{\nb}}{\ra}
F^{p}H^{p+q}_{\dR}(X/\C)
\ot \Omega^r_{\C/\Qb}\\
\os{{\nb}}{\ra}
F^{p-1}H^{p+q}_{\dR}(X/\C)
\ot \Omega^{r+1}_{\C/\Qb},
$$
$$(\text{resp. }
H^{q-1}(\Omega^{p+1}_{X/\C})\ot \Omega^{r-1}_{\C/\Qb}
\os{\ol{\nb}}{\ra}
H^{q}(\Omega^{p}_{X/\C})
\ot \Omega^r_{\C/\Qb}\\
\os{\ol{\nb}}{\ra}
H^{q+1}(\Omega^{p-1}_{X/\C})
\ot \Omega^{r+1}_{\C/\Qb}.
$$

\begin{prop}[\cite{masanori} Proposition 3.4, 3.6]\label{muminv58}
There are the following natural maps:
$$
\Ext^{p}_{\MM(\C)}(\Q(0),H^{q}(X)(r))
\lra
\X_X^{r-p,q-r+p}(p)
\lra
\Lambda_X^{r-p,q-r+p}(p).
$$
\end{prop}
\begin{pf}(sketch)
Note that $\Ext^{p}_{\MM(\C)}(\Q(0),H^{q}(X)(r))=
\lim{S}\Ext^{p}_{\MHM(S)}(\Q_S(0),H^qf_*\Q_{X_S}(r))$.
There is the forgetful functor $\MHM(S)\ra 
{\mathrm{MF}}_{rh}(S)$, which is exact and faithful 
(see \cite{masanori} \S2).
Therefore we have the well-defined map between the Yoneda
extension groups
$$\Ext^{p}_{\MHM(S)}(\Q_S(0),H^qf_*\Q_{X_S}(r))
\ra
\Ext^{p}_{{\mathrm{MF}}_{rh}(S)}(\O_S(0),
R^qf_*\Omega_{X_S/S}^{\bullet}).
$$
Using the Koszul resolution of $\O_{S}$, we can show that the right
hand side is isomorphic to the cohomology at the middle term
of the following complex (cf. \cite{masanori}, Lemma~3.3):
$$
F^{r-p+1}H^{q}_S\ot \Omega^{p-1}_{S/\Qb}
\os{{\nb}}{\ra}
F^{r-p}H^{q}_S
\ot \Omega^p_{S/\Qb}\\
\os{{\nb}}{\ra}
F^{r-p-1}H^{q}_S
\ot \Omega^{p+1}_{S/\Qb},
$$
where we put $F^kH^q_S:=R^qf_*\Omega_{X_S/S}^{\bullet\geq k}$.
Thus we have the natural map
$$
\Ext^{p}_{\MM(\C)}(\Q(0),H^{q}(X)(r))
\lra \lim{S}\Ext^{p}_{{\mathrm{MF}}_{rh}(S)}(\O_S(0),
R^qf_*\Omega_{X_S/S}^{\bullet})
\simeq
\Xi_X^{r-p,q-r+p}(p).$$

The other natural map to $\Lambda_X^{r-p,q-r+p}(p)$
can be constructed similarly.
\end{pf}

\subsection{Filtrations on Chow groups and 
higher Abel-Jacobi maps}\label{higherAJ}
Let $X$ be a nonsingular projective variety over $\C$.
We can construct the cycle map from Bloch's higher Chow
groups (\cite{hc}) to the extension groups in the category of 
arithmetic Hodge modules (\cite{masanori} \S4):
\begin{equation}\label{AHScycle}
\begin{CD}
\lim{X_S}\CH^r(X_S,m)@>>> \lim{X_S}
\Ext^{2r-m}_{\MHM(X_S)}(\Q_{X_S}(0),\Q_{X_S}(r))\\
\Vert@.\Vert\\
\CH^r(X,m)@>>>
\Ext^{2r-m}_{\MM(X)}(\Q_{X}(0),\Q_{X}(r)).
\end{CD}
\end{equation}
We denote the above map \eqref{AHScycle} by $c_X$.

By the Leray spectral sequence \eqref{leray81}, there is the
filtration $F^{\bullet}$ on the right hand side in \eqref{AHScycle}.
We define the filtration on Chow groups as follows:
$$
F^{\nu}\CH^r(X,m):=
c_X^{-1}(F^{\nu}\Ext^{2r-m}_{\MM(X)}(\Q_{X}(0),\Q_{X}(r))).
$$
Since the Leray spectral sequence \eqref{leray81} degenerates
at $E_2$-terms, the cycle map $c_X$ induces the following map.
\begin{equation}\label{hAJ81}
\rho^{\nu}_X:\Gr^{\nu}_F\CH^r(X,m)\lra 
\Ext^{\nu}_{\MM(\C)}(\Q(0),H^{2r-m-\nu}(X)(r)).
\end{equation}
We call the above the {\it $\nu$-th higher Abel-Jacobi map}.
These are injective by definition of $F^{\bullet}\CH^r(X,m)$.

\begin{prop}\label{filt81}
\begin{enumerate}
\renewcommand{\labelenumi}{(\theenumi)}
\item
$F^1\CH^r(X)=\CH^r(X)_{\hom}$. $F^1\CH^r(X,m)=\CH^r(X,m)$ for $m\geq1$.
\item 
$F^2\CH^r(X,m)$ is contained in the kernel of the cycle map to
the Deligne cohomology group.
In particular, we have $F^2\CH^1(X)=0$.
\item
$F^{2}\CH_0(X)=T(X)$.
\item\label{filt81ter}
$F^{r+1}\CH^r(X,m)=F^{r+2}\CH^r(X,m)=\cdots$,
for all $r$.
\item
$F^{\nu}\CH^r(X,m)\cdot F^{\mu}\CH^s(X,n)\subset
F^{\nu+\mu}\CH^{r+s}(X,m+n)$.
\item
$F^{\bullet}$ is respected by any algebraic correspondence.
Each correspondence $\Gamma_*$ on $\Gr_F^{\nu}\CH^r(X,m)$ 
depends only on the K\"{u}nneth $(2\dim X-2r+m+\nu,*)$-component
of the Betti cohomology class $[\Gamma]\in H^*(X\times Y)$.
\end{enumerate}
\end{prop}
\begin{pf}
\noindent(1).
The former follows from that 
$r_{\C}:\MM(\C)\ra \{\Q\text{-vector space}\}$ is
faithful.
The latter follows from the fact that
the cycle class map from $\CH^r(X,m)$ to ordinary Betti cohomology is
zero by a standard weight argument.

\noindent(2). 
It follows from that the realization functor $r_{\C}$ factors through
the category of mixed Hodge structures.

\noindent(3).
Let $\pi=\pi_{2n-1}\in \CH^n(X\times X)$ be the algebraic cycle in 
\cite{murre2} 4.1.Theorem 2, which has the following properties:
(i) $\pi_*^2=\pi_*$, (ii) the cycle class of $\pi$ is the K\"{u}nneth
$(1,2n-1)$ component of the diagonal cycle $\Delta_X$,
(iii)
$\ker(\pi_*:\CH_0(X)_{\dego}\ra \CH_0(X))=T(X)$.
Consider the commutative diagram
$$
\begin{CD}
F^1\CH_0(X) @>{c'}>>
\Ext^{1}_{\MM(\C)}(\Q(0),H^{2n-1}(X)(n))\\
@V{\pi_*}VV@VV{\pi_*}V \\
F^1\CH_0(X) @>>>
\Ext^{1}_{\MM(\C)}(\Q(0),H^{2n-1}(X)(n)).
\end{CD}
$$
Since the property (ii), the right vertical arrow is bijective.
Hence $F^2\CH_0(X)=\ker ~c'$ contains the kernel of the left $\pi_*$,
that is, $T(X)$ by the property (iii).
On the other hand, the Hodge realization functor induces the natural map
$\Ext^{1}_{\MM(\C)}(\Q(0),H^{2n-1}(X)(n))\ra 
\Ext^{1}_{\rMHS}(\Q(0),H^{2n-1}(X)(n))\simeq \Alb(X)(\C)$ which induces
the usual Albanese map, so $F^2\CH_0(X)\subset T(X)$.
Therefore $F^2\CH_0(X)= T(X)$.

\noindent(4).
Let $\nu\geq r+1$ and $n=\dim X$. 
Then there is the following commutative diagram:
$$
\begin{CD}
{\mathrm Gr}^{\nu}_{F}\CH^r(X,m) @>>>
\Ext^{\nu}_{\MM(\C)}(\Q(0),H^{2r-m-\nu}(X)(r))\\
@V{L^{n-2r+m+\nu}}VV@VV{L^{n-2r+m+\nu}}V \\
{\mathrm Gr}^{\nu}_{F}\CH^{n-r+m+\nu}(X,m) @>>>
\Ext^{\nu}_{\MM(\C)}(\Q(0),H^{2n-2r+m+\nu}(X)(n-r)),
\end{CD}
$$
where $L$ is the Lefschetz operator, that is, $x\mapsto x.H$ for a
hyperplane section $H\subset X$.
By the hard Lefschetz theorem, the right vertical arrow is bijective.
On the other hand, the the left vertical arrow is 0 because 
$n-r+m+\nu\geq n+m+1$.
Therefore ${\mathrm Gr}^{\nu}_{F}\CH^r(X)=0$.

\noindent(5).
Clear by definition of $F^{\bullet}\CH^r(X,m)$.

\noindent(6).
It follows from the injectivity of the higher Abel-Jacobi maps 
\eqref{hAJ81}.
\end{pf}

\begin{example}
Let $r=\dim X=n$, $\nu=2$ and $m=0$. Then the higher Abel-Jacobi map gives
the following map:
\begin{equation}\label{alb2}
\rho_X^2:F^2\CH_0(X)=T(X)\lra \Ext^{2}_{\MM(\C)}(\Q(0),H^{2n-2}(X)(n)).
\end{equation}
We call the above \eqref{alb2} the {\it second Albanese map}.
\end{example}

We conjecture that the filtration $F^{\bullet}\CH^r(X,m)$ terminates: 
\begin{conj}\label{mainconj}
$F^{N}\CH^r(X,m)=0$ for some $N\gg 0$. By Proposition 
$\ref{filt81}$ $\eqref{filt81ter}$, it is equivalent to 
$F^{r+1}\CH^r(X,m)=0$. In other words,
the higher Abel-Jacobi map
$$
\rho^{r}_X:F^r\CH^r(X,m)\lra 
\Ext^{r}_{\MM(\C)}(\Q(0),H^{r-m}(X)(r))$$
is injective for each $X$, $r$ and $m$.
\end{conj}
The above conjecture is true for $r=1$.
However, when $r\geq2$, it is a very difficult problem.
\begin{prop}\label{mainbloch326}
If the conjecture \ref{mainconj} is true,  
$F^{\bullet}\CH^r(X,m)$ gives
the conjectural filtration on Chow groups $($cf. \cite{jan},
\cite{murre}, \cite{shuji}$)$. 
In particular, it implies
the Bloch conjecture $\ref{Bloch}$.
\end{prop}
We remark that in order to the Bloch conjecture~\ref{Bloch}, we need only
the injectivity of the second Albanese maps \eqref{alb2} for
any surfaces.
We will show that Conjecture~\ref{mainconj} for $r=2$, $m=0$
follows from the Beilinson conjecture
\ref{Beilinson} for $r=2$ in \S\ref{secBeilinson}.

\medskip

\begin{defn}
Composing the higher Abel-Jacobi map \eqref{hAJ81} and the maps in 
Proposition \ref{muminv58},
we have the following maps:
$$
\xi^{\nu}_X:\Gr^{\nu}_F\CH^r(X,m)
\lra \X_X^{r-\nu,r-m}(\nu),
$$
$$
\delta^{\nu}_X:\Gr^{\nu}_F\CH^r(X,m)
\lra \Lambda_X^{r-\nu,r-m}(\nu).
$$
We call $\xi^{\nu}_X(z)$ (resp. $\delta^{\nu}_X(z)$)
for a cycle $z\in \Gr^{\nu}_F\CH^r(X,m)$, the {\it Mumford $\xi$-invariant
$($resp. $\delta$-invariant$)$ of $z$}.
\end{defn}

For more about the Mumford invariants and its applications,
we are preparing another paper \cite{inf}.

%%%%%%%%%%%%%%%%%%%%%%%%%%%%%%%% Mochizuki %%%%%%%%%%%%%%%%%%%%%%%%%%%%
\section{A calculation of a higher Abel-Jacobi invariant}
The Mumford invariant is not enough to capture all algebraic cycles.
We hope that the invariant obtained from
our higher Abel-Jacobi maps will be stronger and capture all cycles.

In this section, we give an example of a 0-cycle on a surface
whose image by the second Albanese map does not vanish, but so does
the Mumford $\xi$-invariant of it. 
\subsection{}
We use the following curve constructed by Shin-ichi Mochizuki:
\begin{lem}[S.Mochizuki]\label{Mochizuki}
There is a projective nonsingular curve $C$ over $\Qb$ of genus $g\geq 2$,
such that the rank of Neron Severi group of $C\times C$ is $3$, or 
equivalently ${\mathrm End}(J(C))\ot\Q=\Q$.
\end{lem}

We will give the proof of above lemma at the end of this section.

\begin{thm}\label{asakura1}
Let $C$ be a curve as in Lemma~\ref{Mochizuki}.
Let $O\in C(\Qb)$ be a $\Qb$-point of $C$ such that the divisor
$K_C-(2g-2)\cdot O$
is not $\Q$-linearly equivalent to $0$, 
and $P\in C(\C)\setminus\C(\Qb)$ be any $\C$-valued point which is not
contained in the set of $\Qb$-points.
Put $C_{\C}:=C\ot\C$, $X:=C_{\C}\times C_{\C}$ and 
$z:=(P,P)-(P,O)-(O,P)+(O,O)\in T(X)$.
Then we have
$$
\xi_X^{2}(z)=0 \quad \text{in }\quad
\X^{0,2}_X(2),
$$
but
$$
\rho_X^2(z)\not=0 \quad \text{in }
\Ext_{\MM(\C)}^2(\Q(0),H^2(X)(2)).
$$
\end{thm}
\begin{rem}
For any curve $C$ over $\Qb$ of genus $\geq2$, 
there are always $\Qb$-valued points $P$ and
$O$ of $C$ such that the divisor $P-O$ is not $\Q$-linearly equivalent to $0$.
In fact, by the Mumford-Manin conjecture (=Raynaud's theorem), the set 
$C(\Qb)\cap J(C)(\Qb)_{\mathrm tor}$
is finite, which means that the divisor $P-O$ is not torsion except for 
finitely many points $P$.
\end{rem}
\subsection{Proof of Theorem \ref{asakura1}}
First we prove the vanishing $\xi^2_X(z)=0$.
Let $z_S\in \CH^2(S\times C\times C)$ be any model of $z$,
that is, there is a morphism $r:\Spec\C\ra S$ such that
$(r\times 1 \times 1)^*(z_S)=z\in \CH^2(C_{\C}\times C_{\C})$.
By shrinking $S$, we may assume that 
$z_S\vert_{\{t\}\times C \times C}$ is contained in the kernel
of the Albanese map for each $t\in S_{\C}$.
Then, by the following commutative diagram,
$\xi^2_X(z)$ is contained in the image of the map $a_S$
$$
\begin{CD}
T(X) @>{\xi^2_X}>> H^2_{\dR}(X/\C)\ot \Omega^2_{\C/\Qb}
/\nb(F^1H^2_{\dR}(X/\C)\ot \Omega^1_{\C/\Qb})\\
@AAA@AA{a_S}A\\
F^2_S\CH(S\times C\times C) @>>>
R^2p_{S*}\Omega^{\bullet}_{S\times C\times C/S}\ot \Omega^2_{S/\Qb}
/\nb(F^1R^2p_{S*}\Omega^{\bullet}_{S\times C\times C/S}
\ot \Omega^1_{S/\Qb}),
\end{CD}
$$
where $p_S:S\times C\times C\ra C\times C$ denotes the projection.
We can choose the variety $S$ to be $1$-dimensional (see the 
below). Therefore the image of $a_S$ is zero. 
Thus we have $\xi^2_X(z)=0$.

\medskip

Next we show the non-vanishing of the second Albanese class
$\rho^{2}_X(z)$.
By the definition of $\MM(\C)$, the extension group can be written
as an inductive limit of the one in $\MHM(S)$:
$$
\Ext_{\MM(\C)}^2(\Q(0),H^2(X)(2))=
\lim{S}\Ext_{\MHM(S)}^2(\Q_S(0),R^2p_{S*}\Q_{S\times C\times C}(2)).
$$
Therefore $\rho^{2}_X(z)$ does not vanish if and only if so does not
in $\Ext^2_{\MHM(S)}(\Q_S(0),R^2p_{S*}\Q(2))$ 
for any sufficiently dominant $S\ra C$.
The latter condition is equivalent to the non-vanishing in
$\Ext^4_{\MHM(S\times C\times C)}(\Q(0), 
\Q(2))$ by the decomposition theorem of mixed Hodge modules.

Let $f_P:\Spec\C\ra C$ be the associated morphism of the point $P$.
Note that $f_P$ factors through a
generic point of $C$, because $P$ is not a $\Qb$-valued point.
We will take a good lifting $z_0 \in \CH^2(C\times C \times C)$
such that $(f_P\times 1 \times 1)^*(z_0)=z$, and prove that
for any dominant $j_S:S\ra C$, the cycle class of 
$z_S=(j_S\times 1\times 1)^*(z_0)$ does not
vanish in $\Ext^4_{\MHM(S\times C\times C)}(\Q(0),\Q(2))$.

Let $z_0$ be the
{\it modified diagonal cycle} $\Delta_O\in \CH^2(C\times C \times C)$.
Recall the definition of it (\cite{gross}).
Put
$$
\Delta_{xxx}=\{ (x,x,x)\in~C\times C \times C~\vert~x \in C\},
$$
$$
\Delta_{Oxx}=\{ (O,x,x)\in~C\times C \times C~\vert~x \in C\},
$$
$$
\Delta_{xOx}=\{ (x,O,x)\in~C\times C \times C~\vert~x \in C\},
$$
\begin{center}
$\vdots$
\end{center}
the subvarieties of $C\times C \times C$ of codimension 2. 
Then the modified diagonal cycle
$\Delta_O$
is defined as:
$$
\Delta_{xxx}-\Delta_{Oxx}-\Delta_{xOx}-\Delta_{xxO}+\Delta_{OOx}
+\Delta_{OxO}+\Delta_{xOO}.
$$
It is easy to see that $z_0(=\Delta_O)$ is homologically trivial,
and therefore so is $z_S$.
Therefore it defines an Abel-Jacobi class $\rho(z_S)\in
\Ext^1_{\rMHS}(\Q(0),H^1(S_{\C})\ot H^1(C_{\C})^{\ot 2}(2))$.
In order to prove the non-vanishing of the class of $z_S$ in
$\Ext^4_{\MHM(S\times C\times C)}(\Q(0),\Q(2))$, it suffices to show 
the non-vanishing of the class $\rho(z_S)$. 
To do this,
we may assume $\dim ~S=1$ by taking suitable hyperplane
sections.
Thus we have reduced the proof of
Theorem~\ref{asakura1} to the following:
\begin{lem}\label{ken-ichiro}
Let $C_{\C}$, $O$ and $z_0$ as above. Then for any nonsingular
curve $S$ over $\C$ with a dominant morphism $j_S:S\ra C_{\C}$,
the image of the Abel-Jacobi class $\rho(z_0)$ under
the following map
$$(j_S\times 1\times 1)^*:
\Ext^1_{\rMHS}(\Q(0), H^1(C_{\C})^{\ot 3}(2))
\ra
\Ext^1_{\rMHS}(\Q(0),H^1(S)\ot H^1(C_{\C})^{\ot 2}(2))
$$
does not vanish.
\end{lem}
\begin{pf}
We write $C_{\C}$ by $C$ simply.
We prove the assertion in the following step:
\begin{enumerate}
\renewcommand{\labelenumi}{(\ref{ken-ichiro}.\theenumi)}
\item\label{aaa}
The Abel-Jacobi class 
$\rho(z_0)\in \Ext^1_{\rMHS}(\Q(0), H^1(C_{\C})^{\ot 3}(2))$
is not zero.
\item\label{bbb}
The image of the Abel-Jacobi class $\rho(z_0)$ under
the following map
$$
\Ext^1_{\rMHS}(\Q(0), H^1(C)^{\ot 3}(2))
\ra
\Ext^1_{\rMHS}(\Q(0),H^1(C)\ot \Sym^2(H^1(C))(2))
$$
does not vanish.
Here $\Sym^{2}(H)=H^{\ot2}/\{a\ot b-b \ot a\}$ denotes
the symmetric product.
\item\label{ccc}
The natural map
$$
\Ext^1_{\rMHS}(\Q(0),H^1(C)\ot \Sym^2(H^1(C))(2))
\ra
\Ext^1_{\rMHS}(\Q(0),H^1(S)\ot \Sym^2(H^1(C))(2))
$$
is injective.
\end{enumerate}

\bigskip

\noindent{\it Proof of} (\ref{ken-ichiro}.\ref{aaa}).
Let $p_i:C\times C\times C\ra C\times C$ be
the $i$-th projection,
and
$p_{ij}:C\times C\times C\ra C\times C$ 
the projection into $(i,j)$-th component.
We note that $\Delta_{xxx}=p_{12}^*\Delta.p_{23}^*\Delta$,
$\Delta_{Oxx}=p_{23}^*\Delta.p_{1}^*O$, $\cdots$,
where $\Delta\subset C\times C$ the diagonal cycle.

Let $i:C\times C \hra C\times C \times C$ be an inclusion such as:
$(x,y)\mapsto (x,y,y)$.
We can easily show that the pull-back of the modified diagonal cycle 
$i^*z_0$ is equal to $\Delta^2+O\times K_C$ in $\CH_0(C\times C)$.
We want to show the non-vanishing of $i^*(\rho(z_0))$. It is the image
of $\rho(z_0)$ under the map 
$$
\Ext^1_{\rMHS}(\Q(0), H^1(C)^{\ot 3}(2))
\lra
\Ext^1_{\rMHS}(\Q(0),H^1(C)\ot H^2(C)(2))
$$
induced from the cup-product $H^1(C)^{\ot2}\ra H^2(C)$. 
Under the isomorphism $H^1(C)\ot H^2(C)(2)\simeq H^1(C)(1)$, 
$i^*(\rho(z_0))$ coincides with the Abel-Jacobi class of 
$q_{1*}i^*z_0=-K_C+(2g-2)\cdot O\in \CH_0(C)$ 
where $q_1:C\times C \ra C$ denotes the 1st projection.
Since $K_C-(2g-2)\cdot O$ is not $\Q$-linearly
equivalent to $0$,
it does not vanish.

\medskip

\noindent{\it Proof of} (\ref{ken-ichiro}.\ref{bbb}).
More generally, we show the image of $\rho(z_0)$ under the map
$$
\Ext^1_{\rMHS}(\Q(0), H^1(C)^{\ot 3}(2))
\ra
\Ext^1_{\rMHS}(\Q(0),\Sym^3(H^1(C))(2))
$$
does not vanish. Put $\ol{\rho(z_0)}$ be the image.
Let $\alpha:\Sym^3(H^1(C))\ra H^1(C)^{\ot 3}$ be
a morphism of Hodge structures defined as: $v_1\cdot v_2\cdot v_3
\mapsto \sum_{\sigma\in {\frak S}_3}v_{\sigma(1)}\ot v_{\sigma(2)}
\ot v_{\sigma(3)}$.
Then we want to show $\alpha(\ol{\rho(z_0)})\not=0$.
Since the modified diagonal cycle is invariant under the action 
of the symmetric group ${\frak S}_3$, we have 
$\alpha(\ol{\rho(z_0)})=\sum_{\sigma\in {\frak S}_3}\rho(z_0)^{\sigma}
=\sum_{\sigma\in {\frak S}_3}\rho(z_0^{\sigma})
=6\rho(z_0)$. 
Therefore the assertion follows from (\ref{ken-ichiro}.\ref{aaa}).

\medskip

\noindent{\it Proof of} (\ref{ken-ichiro}.\ref{ccc}).
Let $\ol{S}$ be a smooth completion of $S$, and $\ol{j}_S:\ol{S}\ra C$
be the extension of the morphism $j_S$.

Since the category of polarizable pure Hodge structures is semi-simple,
the map
$$
\Ext^1_{\rMHS}(\Q(0),H^1(C)\ot \Sym^2(H^1(C))(2))
\ra%\os{(\ol{j}_S\times 1 \times 1)^*}{\lra}
\Ext^1_{\rMHS}(\Q(0),H^1(\ol{S})\ot \Sym^2(H^1(C))(2))
$$
is injective. Therefore it suffices to show that
$$
\Ext^1_{\rMHS}(\Q(0),H^1(\ol{S})\ot \Sym^2(H^1(C))(2))
\lra
\Ext^1_{\rMHS}(\Q(0),H^1(S)\ot \Sym^2(H^1(C))(2))
$$
is injective.
There is the exact sequence of mixed Hodge structures:
$$
0\lra H^1(\ol{S},\Q)\lra H^1(S,\Q)\lra \op\Q(-1)\lra 0.
$$
Therefore,
applying ${\Bbb R}\Hom_{\rMHS}(\Q(0),(-)\ot \Sym^2(H^1(C))(2))$ to this
short exact sequence, our assertion follows from that
$\Hom_{\mathrm HS}(\Q(0),\Sym^2(H^1(C))(1))=0$.
Consider the surjection:
$$
\Hom_{\mathrm HS}(\Q(0),(H^1(C))^{\ot2}(1))\lra
\Hom_{\mathrm HS}(\Q(0),\Sym^2(H^1(C))(1)).
$$
The left hand side is isomorphic to ${\mathrm End}(J(C))\ot\Q$, 
which is a
1-dimensional vector space generated
by the diagonal cycle class $[\Delta]$ because of
$\rank \NS(C\times C)=3$.
The class $[\Delta]$ vanishes into the right hand side,
which means $\Hom_{\mathrm HS}(\Q(0),\Sym^2(H^1(C))(1))=0$.

\medskip

Thus we have proved all steps.
\end{pf}
\begin{rem}[A `` heuristic proof " of Theorem~\ref{asakura1}]
The most technical 
point of the above proof is to choose the modified diagonal cycle
as a lifting of the 0-cycle
$z$.
The reader may have a question why we choosed it.
We give a heuristic answer to this question.

\medskip

Let $z_0\in \CH^2(C\times C\times C)$ be a lifting of $z$ which is 
homologically trivial. Then we have the Abel-Jacobi class of 
$\rho(z_0)\in \Ext^1_{\rMHS}(\Q(0),H^1(C_{\C})^{\ot3}(2))$.
The essential part of the above proof is Lemma \ref{ken-ichiro}, that is,
to show 
that the image of $\rho(z_0)$ under the natural map
$$
\Ext^1_{\rMHS}(\Q(0), H^1(C_{\C})^{\ot 3}(2))
\ra
\Ext^1_{\rMHS}(\Q(0),H^1(S)\ot H^1(C_{\C})^{\ot 2}(2))
$$
does not vanish for all dominant maps $S\ra C_{\C}$ from
nonsingular curves $S$.

Consider the following commutative diagram:
$$
\begin{CD}
\Ext^1_{\rMHS}(\Q(0), H^1(C_{\C})^{\ot 3}(2))
@>>>
\Ext^1_{\rMHS}(\Q(0),H^1(S)\ot H^1(C_{\C})^{\ot 2}(2))
\\
@V{\alpha}VV @VVV \\
\Ext^1_{\rMHS}(\Q(0), H^1(C_{\C})^{\ot 3}/N^1(2))
@>{\beta}>>
\Ext^1_{\rMHS}(\Q(0),H^1(\ol{S})\ot H^1(C_{\C})^{\ot 2}/N^1(2)),
\end{CD}
$$
where $N^{\bullet}$ denotes the coniveau filtration (cf. \cite{jan} p.162).
The map $\beta$ is injective.
Therefore we only have to choose a lifting $z_0$ satisfying 
$\alpha(\rho(z_0))\not=0$. 
This is a well-known problem concerning with the {\it Griffiths groups}.
The Griffiths group of a projective nonsingular variety $X$ (denoted by
${\mathrm Grif}^r(X)$) is defined to be the group of homologically trivial
cycle modulo algebraic equivalence: ${\mathrm Grif}^r(X):=\CH^r(X)_{\hom}/
\CH^r(X)_{\mathrm alg}$. In the conjectural theory of mixed motives, 
there is the following isomorphism:
$$
{\mathrm Grif}^r(X)\simeq
\Ext^1_{\mot(\C)}(\Q(0), H^{2r-1}(X)/N^{r-1}(r)).
$$
The right hand side is conjectured to be injected into the extension group
of mixed Hodge structures
$\Ext^1_{\rMHS}(\Q(0), H^{2r-1}(X)/N^{r-1}(2))$.
Therefore to find a cycle $z_0$ such that $\alpha(\rho(z_0))\not=0$ is 
conjecturally equivalent to find a non-trivial element 
in the Griffiths group.
I was taught from K.Kimura that the modified diagonal cycle for certain curves 
gives the non-trivial element of the Griffiths groups (cf. \cite{kimura}).
But we do not necessarily need 
to prove it for Lemma \ref{ken-ichiro}.
I don't know whether the modified diagonal cycle gives a non-trivial
element of the Griffiths group for Mochizuki's curve 
(Lemma~\ref{Mochizuki}).
\end{rem}
\begin{rem}
When $C$ is an elliptic curve, the 0-cycle $z=(P,P)-(P,O)-(O,P)+(O,O)$ on
$C\times C$ is rationally equivalent to 0. (Easy exercise. Hint: 
consider the embedding $C\hra C\times C$, $x\mapsto (x,-x+P)$.)
\end{rem}
\subsection{Proof of Lemma \ref{Mochizuki}}\label{m55}
We can easily find the desired curve at least over $\C$:
\begin{lem}\label{m1}
There is a nonsingular projective curve $X$ of genus $g\geq2$ 
defined over
$\C$ such that $\End(J(X))=\Z$.
\end{lem}
\begin{pf}
Let $S_0$ be a nonsingular projective surface over $\Qb$ 
with irregularity $q=0$
(e.g. $S_0=\P^2_{\Qb}$), and $L$ be a very ample line bundle on $S_0$.
We assume that
general smooth member of the linear system $\vert L\vert$ is not
hyper-elliptic.
We consider a Lefschetz pencil $f:S \ra \P^1_{\Qb}$ obtained from $L$.
We will show that the 
geometric generic fiber $C_{\ol{\eta}}=f^{-1}(\ol{\eta})$ 
satisfies $\End(J(C_{\ol{\eta}}))=\Z$.

Note that $\End(J(C_{\ol{\eta}}))\subset\End
(H^1_{\text{\'{e}t}}(C_{\ol{\eta}},\Q_{\l}))$, and $\gamma\in 
\End(J(C_{\ol{\eta}}))$ satisfies
$\gamma\cdot T_j^{n_j}= T_j^{n_j}\cdot\gamma$ for some $n_j>0$, where
$T_j\in \pi_1(\eta,\ol{\eta})$ is the local monodromy generator.
Let $\delta_j$ be the 
corresponding vanishing cycle.
Then by the Picard-Lefschetz formula, we have
$$
T_j(x)=x-(x,\delta_j)\delta_j \quad  \text{ for } x  
\in H^1_{\text{\'{e}t}}(C_{\ol{\eta}}).
$$
Therefore we have
$(x,\delta_j)\gamma(\delta_j)=(\gamma(x),\delta_j)\delta_j$ for all 
$x \in H^1_{\text{\'{e}t}}(C_{\ol{\eta}})$. We put
\begin{equation}\label{van}
\gamma(\delta_j)=\ve_j\delta_j
\end{equation}
for some $\ve_j \in \Q_{\l}$.
Moreover, applying $T_i^{n_i}$ to \eqref{van}, we have
$(\delta_i,\delta_j)(\ve_i-\ve_j)=0$ for all $i$, $j$.
Since every vanishing cycles are conjugate under the action of 
$\pi_1(\eta,\ol{\eta})$,
(that is, for any $\delta_i$ and $\delta_j$, 
there is a $\sigma \in \pi_1(\eta,\ol{\eta})$ such
that $\sigma(\delta_i)=a\cdot\delta_j$ for some $a \in\Q^*_{\l}$), 
we can see that
for any $\delta_i$ and $\delta_j$, there are $\delta_{k_1},\cdots,\delta_{k_l}$
such that $(\delta_{i},\delta_{k_1})\not=0,
(\delta_{k_1},\delta_{k_2})\not=0,\cdots, (\delta_{k_l},\delta_{j})\not=0$.
Therefore all the $\ve_i$ coincide. 

The vanishing cycles generate 
$H^1_{\text{\'{e}t}}(C_{\ol{\eta}})$ because of $q=0$.
Therefore we have $\gamma=n\cdot \mathrm{id}$, which implies the assertion.
\end{pf}
\def\M{{\cal M}}
\def\Mb{\ol{\cal M}}

\begin{lem}\label{m4}
There is an absolutely unramified $p$-adic number field $K$,
and a proper curve 
\begin{equation}\label{serre51}
\pi:X\lra \Spec \O_K
\end{equation}
over the integer ring $\O_K$ such that,
\begin{enumerate}
\renewcommand{\labelenumi}{(\theenumi)}
\item\label{1-83}
$\pi$ is smooth over $K$, and $\End(J(X_{\ol{K}}))=\Z$,
\item\label{2-83}
The special fiber $X_0$ is a geometrically 
connected stable curve $E_1\cup\cdots\cup E_g$ over $\F_q$ $(=$ 
the residue field, $q=p^n$ $)$ such that
\begin{enumerate}
\renewcommand{\labelenumi}{(\alph{\enumi})}
\item\label{2a-83}
each $E_i$ is an ordinary elliptic curve, that is,
$\End(E_i)\ot\Q$ is a quadratic number field, and
the rational prime $p$ is complete split $:$ $p=\fp_1\fp_2$,
$($Then $\End(E_i)=\End(E_{i,\ol{\F}_q})$ $)$
\item\label{2b-83}
$E_i$ and $E_{i+1}$ intersect at one point for $i=1,\cdots, g-1$,
and
$E_i$ and $E_{j}$ do not intersect if $\vert j-i\vert \geq 2$,
\item\label{2c-83}
$\Hom(E_i,E_j)=0$ for $i\not=j$.
\end{enumerate}
\end{enumerate}
\end{lem}
\begin{pf}
Let $E_{1,F},\cdots,E_{g,F}$ be geometrically connected CM-elliptic curves
over a number field $F$ such that
\begin{enumerate}
\renewcommand{\theenumi}{\roman{enumi}}
\renewcommand{\labelenumi}{(\theenumi)}
\item\label{1-85}
$\Hom(E_{i,F},E_{j,F})=0$ for $i\not=j$,
\item\label{2-85}
$\End(E_{i,F})=\End(E_{i,\Qb})$, 
\item\label{3-85}
there is the smooth model $E_{i,\O_F}\ra \Spec\O_F$ over
the integer ring $\O_F$ of $F$.
\end{enumerate}
Let $Y_{\O_F}=E_1^{\O_F}\cup\cdots\cup E_g^{\O_F}$ be a chain of the
above elliptic curves such as in the condition \eqref{2b-83}.

Let $\M_g$ be the moduli scheme of 
nonsingular proper curves of genus $g\geq2$
over $\O_F$, and $\Mb_g$ the compactification of it in the sense of
\cite{DM}.
Note that $\Mb_g$ is not a scheme but a stack.
There is an integer $N\not=0$, an \et open $U\ra \Mb_g$ where $U$ is a 
regular scheme, and a point $y:\Spec\O:=\Spec\O_F[1/N]\ra U$
which associates with $Y_{\O}:=Y_{\O_F}\ot_{\O_F}\O_F[1/N]$.
We assume that $\Spec\O_F[1/N]\ra \Spec\Z$ is \et by replacing $N$ by
a suitable one.

\smallskip

We denote $F_{\fp}$ the completion of $F$ by a prime $\fp$, and
$\O_{F_{\fp}}$ its integer ring.

\begin{claim}
Let $\fp\in\Spec\O$ be any prime, and $\O_{y(\fp)}=\O_{U,y(\fp)}$ the stalk
of $\O_U$ at the point $y(\fp)$.
Then,
there is an embedding
$$
\alpha:\O_{y(\fp)}\hra \O_{F_{\fp}}
$$
of $\O_F$-local rings 
which induces the isomorphism of the residue fields.
\end{claim}
\begin{pf}
Let $\hat{\O}_{y(\fp)}$ be the completion of $\O_{y(\fp)}$ by
the maximal ideal $\fm$.
Since $\O_{y(\fp)}$ is a regular local ring, there is an
isomorphism $\hat{\O}_{y(\fp)}\simeq \O_{F_{\fp}}[[t_1,\cdots,t_{3g-3}]]$.
We may assume that each $t_j$ is contained in $\O_{y(\fp)}$.
In fact, we can replace $t_j$ by $t'_j$ such that $t_j-t'_j\in \fm^2$.

Since the transcendental degree of $F_{\fp}$ over $F$ is infinite,
there are $(3g-3)$-elements $q_1,\cdots,q_{3g-3}\in \fp\O_{F_{\fp}}$
which are algebraically independent.
Then the $\O_{F_{\fp}}$-algebra homomorphism
$\hat{\O}_{y(\fp)}=
\O_{F_{\fp}}[[t_1,\cdots,t_{3g-3}]]\ra \O_{F_{\fp}}$ defined by 
$t_j\mapsto q_j$ induces the $\O_F$-algebra homomorphism
$\alpha:\O_{y(\fp)}\ra \O_{F_{\fp}}$. 

The remaining part of the proof is to show that $\alpha$ is injective.
Let $I$ be the kernel of $\alpha$.
Since $\O_{y(\fp)}\ot_{\O_F}F\ra \O_{F_{\fp}}\ot_{\O_F}F$ is a homomorphism
of fields (and hence injective), we have $I\ot_{\O_F}F=0$.
On the other hand, $I\subset \O_{y(\fp)}\subset 
\O_{F_{\fp}}[[t_1,\cdots,t_{3g-3}]]$ is torsion free over $\O_F$.
Therefore we have $I=0$.
\end{pf}

Let $$X_{\O_{F_{\fp}}}\ra\Spec\O_{F_{\fp}}\quad(\fp\in\Spec\O)$$
be the proper curve over $\O_{F_{\fp}}$ which is the pull-back of
the family $X_U\ra U$ by $\alpha$.
By Lemma~\ref{m1}, the generic fiber satisfies the condition
\eqref{1-83}.
The special fiber satisfies \eqref{2b-83} and \eqref{2c-83}.
Moreover, there are infinite many prime $\fp$ such that
each irreducible component of the special fiber is ordinary.
In fact, the density of such primes is larger than or equal to 
$1/2^g$.
So the condition \eqref{2a-83} is also satisfied.

This completes the proof.
\end{pf}

\begin{lem}\label{m5}
Let $A\ra \Spec R$ be a proper smooth abelian scheme over a discrete
valuation ring $R$. Put $A_n:=A\ot_RR/\fm^{n+1}$ where $\fm$ is the maximal
ideal.
Then the natural maps $\End(A)\ra \End(A_0)$
and $\End(A_n)\ra \End(A_0)$ are injective.
\end{lem}
\begin{pf}
Due to $$\End(A)\subset\End(A\ot_R\hat{R})=\us{n}{\varprojlim}~\End(A_n)
\quad
(\hat{R}:=\us{n}{\varprojlim}R/\fm^{n+1}),$$
it suffices to show the latter assertion.

Let $\l$ be a rational prime which is invertible in $R/\fm$.
Then the group scheme ${}_{\l^r}(A_n)$ of the $\l^r$-torsion points
is finite \et over $\Spec R/\fm^{n+1}$, and therefore so is the endomorphism
scheme ${\cal{E}}\text{\it{nd}}({}_{\l^r}(A_n))$.
By the formal \'etaleness (\cite{EGA} \S17), we have
\begin{multline*}
\End({}_{\l^r}(A_n))=\Hom(\Spec R/\fm^{n+1},
{\cal{E}}\text{\it{nd}}({}_{\l^r}(A_n)))\\
\simeq \Hom(\Spec R/\fm,{\cal{E}}\text{\it{nd}}({}_{\l^r}(A_n)))=\End({}_{\l^r}(A_0)).
\end{multline*}
Therefore it suffices to show that the natural map
$\End((A_n))\ra{\prod}_r\End({}_{\l^r}(A_n))$ is injective.
It follows from the fact that the subgroup of the $\l$-primary torsion
points is schematically dense.
\end{pf}
Let the notations as in Lemma~\ref{m4}.
The Jacobian variety $J(X_0)$ of $X_0$ is isomorphic to
$E_1\times\cdots\times E_g$, and $\End(J(X_0))$
is isomorphic to
$\End(E_1)\times\cdots\times \End(E_g)$.

\begin{lem}\label{m2}
Consider an arbitrary lifting $X'\ra \Spec \O_K$ of $X_0$.
Let $\varphi$ be an endomorphism of $J(X_0)$. If there is a
positive integer $m>0$ such that $m\cdot\varphi$ can be lifted
on an endomorphism of $J(X'_{\ol{K}})$, then so does $\varphi$:
$$\End(J(X'_{\ol{K}}))=(\End(J(X'_{\ol{K}}))\ot\Q)\cap 
\End(J(X_{0})).$$
\end{lem}
\begin{pf}
Let $\psi:=\wt{m\cdot \varphi}\in \End(J(X'_{\ol{K}}))$ 
be the lifting of $m\cdot \varphi$.
We want to show that
\begin{equation}\label{tate83}
\psi^*(H^1_{\text{\'{e}t}}(J(X'_{\ol{K}}),\Z_{\l})
\subset m\cdot H^1_{\text{\'{e}t}}(J(X'_{\ol{K}}),\Z_{\l})
\end{equation}
for all rational prime $\l$.
If $\l\not=p$, it follows from the isomorphism 
$H^1_{\text{\'{e}t}}(J(X'_{\ol{K}}),\Z_{\l})\simeq
H^1_{\text{\'{e}t}}(J(X_{0,\ol{\F}_q}),\Z_{\l})$. 
So we may assume $\l=p$.

Firstly, we note that $\psi$ is defined
over $K$.
In fact, assume that $\psi$ is defined over a finite extension
$K'$ over $K$. We may assume that $K'/K$ is a Galois extension.
Then,
for any $\sigma\in \Gal(K'/K)$, $\psi^{\sigma}$ is also a lifting
of $m\cdot \varphi$, which coincides with $\psi$ by Lemma~\ref{m5}.

Since the crystalline cohomology depends only on 
the special fiber, we have 
$$\psi^*(H^1_{\crys}(J(X_{0})/\O_K))
\subset m\cdot H^1_{\crys}(J(X_{0})/\O_K).$$
Moreover, by the isomorphism
$H^1_{\crys}(J(X_{0})/\O_K)
\simeq
H^1_{\dR}(J(X')/\O_K)
$ (because $K$ is absolutely unramified),
we have 
$$\psi^*(H^1_{\dR}(J(X')/\O_K))
\subset m\cdot H^1_{\dR}(J(X')/\O_K).$$
Therefore the assertion \eqref{tate83} follows by applying
the Dieudonne functor (\cite{faltings}, \cite{fontaine} 9.11).  
\end{pf}

\begin{lem}\label{m3}
There are finitely many subalgebras $A_0,\cdots,A_m$ of $\End(J(X_0))$
such that, for any lifting $X'\ra \Spec \O_K$ of $X_0$, the image of 
$\End(J(X'_{\ol{K}}))$ is $A_i$ for some $i$.
\end{lem}
\begin{pf}
By Lemma~\ref{m2}, we have
$$\End(J(X'_{\ol{K}}))=(\End(J(X'_{\ol{K}}))\ot\Q)\cap \End(J(X_0)),$$
where $\End(J(X'_{\ol{K}}))\ot\Q$ is a $\Q$-subalgebra of 
$\End(J(X_{0,\ol{\F}_q}))\ot\Q=\End(J(X_0))\ot\Q$.
Since $\End(J(X_0))\ot\Q$ is a reduced and finite-dimensional $\Q$-algebra,
there are at most finitely many $\Q$-subalgebras of it. 
Thus we have the assertion.
\end{pf}

Let us prove Lemma~\ref{Mochizuki}.

\vspace{0.1cm}

Consider a deformation $X'\ra \Spec\O_K$ for each integer
$m\geq0$ such that
\begin{itemize}
\item
$X' \ot_{\O_K}\O_K/p^{m+1}=X\ot_{\O_K}\O_K/p^{m+1}$,
\item
the generic fiber is defined over a number field.
\end{itemize}
By Lemma~\ref{m3}, $\End(J(X'_{\ol{K}}))=A_i$ for some $i\geq0$.
We may assume that $A_0=\Z$.
We want to show that 
$\End(J(X'_{\ol{K}}))=A_0=\Z$ if $m$ is suffciently large.

There is a finite extension $K'/K$ such that 
$\End(J(X'_{K'}))=\End(J(X'_{\ol{K}}))$.
We put $X_n:=X\ot_{\O_{K}}\O_{K'}/p^{n+1}$ and 
$X'_n:=X'\ot_{\O_{K}}\O_{K'}/p^{n+1}$.
Choose an endomorphism $\varphi_i\in A_i\setminus\Z$ for each $i\geq1$.
Since $\End(J(X_{K'}))= \us{n}{\cap}\End(J(X_n))$, there is an integer
$N>0$ such that $\varphi_i\not\in\End(J(X_{N}))$ for all $i\geq1$.
Let $m\geq N$.  
Since
$\End(J(X'_{K'}))\subset \End(J(X'_{N}))=\End(J(X_{N}))$, 
we have $\varphi_i\not\in\End(J(X'_{K'}))$, which
means $\End(J(X'_{\ol{K}}))=A_0=\Z$.

This completes the proof of Lemma~\ref{Mochizuki}.
%%%%%%%%%%%%%%%%%%%%%%%%%%%%%%%%%% Beilinson %%%%%%%%%%%%%%%%%%%%%
\section{The Beilinson conjecture}\label{secBeilinson}
\subsection{}
Recall the {\it Beilinson conjecture} (\cite{jan} 11.4, c)):
\begin{conj}[Beilinson]\label{Beilinson}
Let $r\geq2$ be an integer.
\begin{enumerate}
\renewcommand{\labelenumi}{(\theenumi)}
\item\label{Bei310}
The Abel-Jacobi map
\begin{equation}
\rho:\CH^r(X)_{\hom}\lra J^r(X_{\C})
\end{equation}
is injective for any nonsingular projective variety $X$ over $\Qb$.
\item\label{comparison}
Let $X$ be a nonsingular projective variety $X$ over $\Qb$, and
$z \in \CH^r(X_{\C})$ be any algebraic cycle on $X_{\C}=X\ot_{\Qb}\C$. 
Then, for each $\sigma\in \Aut(\C/\Qb)$, 
$\rho(z)=0$ if and only if $\rho(z^{\sigma})=0$.
\end{enumerate}
\end{conj}
For a projective nonsingular variety $X$ over $\C$,
it is conjectured that the motivic filtration $F^2_{\cal M}\CH^r(X)$
coincides with the kernel of the Abel-Jacobi map
$\rho:\CH^r(X)_{\hom}\ra J^r(X)$.
The Beilinson conjecture~\ref{Beilinson}~(1) 
asserts that $F_{\cal M}^2\CH^r(X)=0$
if $X$ is defined over a number field.
The Beilinson conjecture~\ref{Beilinson}~(2) asserts that
the kernel of Abel-Jacobi maps should be ``algebraic''.
It holds, at least, when $r=1$ or $\dim X$. (I do not
know whether it holds in case $r\not=1,\dim X$.)

\begin{thm}\label{asakura2}
The Beilinson conjecture~$\ref{Beilinson}$~$(1)$, $(2)$ for $r=2$
implies the Bloch conjecture
$\ref{Bloch}$.
\end{thm}
\begin{pf}
By Proposition~\ref{mainbloch326} and the remark after it, 
we show Conjecture~\ref{mainconj} for $r=2$, $m=0$
under the Beilinson conjecture \ref{Beilinson} (1), (2) for $r=2$.
By the definition of the cycle map \eqref{AHScycle}, it suffices to show
that the following map
$$
\CH^2(X_S)\lra \Ext^4_{\MHM(X_S)}(\Q_{X_S}(0),\Q_{X_S}(2))
$$
is injective for any model $X_S$.
If $X_S$ is projective, it follows directly from the Beilinson conjecture 
\ref{Beilinson} (1).
For an open case, we will use the Beilinson conjecture 
\ref{Beilinson} (2).

%\medskip

Let $U$ be a nonsingular quasi-projective (not necessarily complete) 
variety over $\Qb$,
and $Y$ be a smooth completion of $U$ such that $D=\cup D_i=Y-U$
is a simple normal crossing divisor.
We want to show that the map 
$$
\CH^2(U)\lra \Ext^4_{\MHM(U_{\C})}(\Q(0),\Q(2))
$$
is injective.
Since any extension groups of degree$\geq2$
in $\rMHS$
vanishes, the right hand side is an extension of $W_4H^4(U_{\C},\Q)\cap
F^2$ by
$\Ext_{\rMHS}^1(\Q(0),H^3(U_{\C},\Q(2))$.
Put $\CH^2(U)_{\hom}:=\ker (\CH^2(U)\ra H^4(U_{\C},\Q))$.
We will show that
\begin{equation}\label{inj34}
\CH^2(U)_{\hom}\lra\Ext_{\rMHS}^1(\Q(0),H^3(U_{\C},\Q(2))
\end{equation}
is injective. 

Let $\CH_Y^1(D)_{\hom}$ be the subgroup of $\us{i}{\op}\CH^1(D_i)$ 
generated by cycles which are homologicaly equivalent to $0$ in $Y$, 
that is, 
$\CH^1_Y(D)_{\hom}
=\ker(\us{i}{\op}\CH^1(D_{i})\ra
\CH^2(Y)/\CH^2(Y)_{\hom})$.
Consider the following commutative diagram:
\begin{equation}
\begin{CD}
@.\text{(exact)}\\
@.0\\
@.@VVV\\
\CH^1_Y(D)_{\hom}@>{\rho_{D}}>>
J(Y_{\C},D_{\C}) \\
@V{i}VV@VV{a}V\\
\CH^2(Y)_{\hom}@>{\rho_Y}>>\Ext_{\rMHS}^1(\Q(0),H^3(Y_{\C},\Q(2))\\
@VVV@VV{b}V\\
\CH^2(U)_{\hom}@>{\rho_U}>> \Ext_{\rMHS}^1(\Q(0),H^3(U_{\C},\Q(2)))\\
@VVV\\
0@.\\
\text{(exact).}
\end{CD}
\end{equation}
Here we put $J(Y_{\C},D_{\C})=\ker b$.
Note that $\rho_Y$ is injective by the Beilinson 
conjecture~\ref{Beilinson}~\eqref{Bei310}.
We show that if an algebraic cycle $z\in \CH^2(Y)_{\hom}$ satisfies
$b(\rho_Y(z))=0$, then there is a cycle $w\in\CH^1_Y(D)_{\hom}$
such that $\rho_Y(z)=a(\rho_{D}(w))=\rho_Y(i(w))$.
\begin{lem}\label{kernel323}
$J(Y_{\C},D_{\C})=\rho_{D_{\C}}(\CH^1_{Y_{\C}}({D_{\C}})_{\hom})$, where
$\CH^1_{Y_{\C}}({D_{\C}})_{\hom}$ 
denotes the subgroup of $\us{i}{\op}\CH^1(D_{i,\C})$ 
generated by cycles which are homologicaly equivalent to $0$ in $Y_{\C}$
\end{lem}
\begin{pf}
\def\J#1{J({#1})}
For a $\Q$-mixed Hodge structure $H$,
we write $\Ext_{\rMHS}^1(\Q(0),H)$ by $\J H$ simply.

\def\D{{D_{\C}}}
\def\Di{{D_{i,\C}}}
\def\Y{{Y_{\C}}}
\def\U{{U_{\C}}}

There are the following exact sequences:
\begin{equation}\label{exact324}
\Gr^W_3H^3_{\D}(\Y)\lra H^3(\Y)\lra W_3H^3(\U)\lra 0,
\end{equation}
\begin{equation}\label{exact2324}
0\lra W_3H^3(\U)\lra W_4H^3(\U)\lra \Gr^W_4H^3(\U)\lra 0.
\end{equation}
(Here $H^{\bullet}(X_{\C})(r)$ denotes the Betti cohomology 
$H^{\bullet}(X_{\C}^{an},\Q(r))$.)
By an easy computation, we can see that there is a natural isomorphism
$\Gr^W_3H^3_{\D}(\Y)\simeq \us{i}{\op}H^1(\Di)(-1)$ of Hodge structure, 
and a surjection
$\ker(\us{i}{\op}H^2(\Di)(-1)\ra H^4(\Y))\ra \Gr^W_4H^3(\U)$.
We put $V=\ker(\us{i}{\op}H^2(\Di)(-1)\ra H^4(\Y))$.
We have the exact sequence
\begin{equation}\label{Jexact324}
\us{i}{\op}\J{H^1(\Di)(1)}\lra \J{H^3(\Y)(2)}\lra \J{W_3H^3(\U)(2)}\lra 0
\end{equation}
from \eqref{exact324}, and
\begin{equation}\label{Jexact2324}
\Hom_{\mathrm HS}(\Q(0),V(2))
\lra \J{W_3H^3(\U)(2)}\lra \J{W_4H^3(\U)(2)}(=\J{H^3(\U)(2)})
\end{equation}
from \eqref{exact2324}.
By the above two exact sequences \eqref{Jexact324} and \eqref{Jexact2324}, 
we have the following exact sequence
\begin{equation}\label{JYexact324}
\us{i}{\op}\J{H^1(\Di)(1)}\lra J(\Y,\D)\lra 
\Hom_{\mathrm HS}(\Q(0),V(2))\lra 0.
\end{equation}
By the Lefschetz $(1,1)$ theorem, 
$\Hom_{\mathrm HS}(\Q(0),V(2))=\ker(\us{i}{\op}\NS(\Di)\ra H^4(\Y))$.
We thus have the following commutative diagram:
$$
\begin{CD}
0\\
@VVV\\
\us{i}{\op}\CH^1(\Di)_{\hom}@>{\op \rho_{\Di}}>>\us{i}{\op}\J{H^1(\Di)(1)}\\
@VVV@VVV\\
\CH^1_{\Y}(\D)_{\hom}@>{\rho_{\D}}>>J(\Y,\D)\\
@VVV@VVV\\
\ker(\us{i}{\op}\NS(\Di)\ra H^4(\Y))@>{\simeq}>>\Hom_{\mathrm HS}(\Q(0),V(2))\\
@VVV@VVV\\
0 @.0.
\end{CD}
$$
Since the Abel-Jacobi map $\rho_{\Di}$ is bijective,
we obtain the surjectivity of $\rho_{\D}$.
\end{pf}

\medskip

By Lemma~\ref{kernel323}, there is an algebraic cycle 
$w_{\C}\in \CH^1_{Y_{\C}}(D_{\C})_{\hom}$ 
such that $\rho_{Y_{\C}}(z)=a(\rho_{D_{\C}}(w_{\C}))=\rho_{Y_{\C}}(i(w_{\C}))$.

The remaining part of the proof is to show that there is a cycle 
$w\in\CH^2(Y)_{\hom}$ (defined over $\Qb$) with support on $D$ such that
$\rho_{Y_{\C}}(i(w_{\C}))=\rho_{Y_{\C}}(w)$.
We write $i(w_{\C})$ by $w_{\C}$ simply.
To do this, we use the Beilinson conjecture \ref{Beilinson} (2).
Since the cycle $z-w_{\C}$ is contained in the kernel of the Abel-Jacobi
map $\rho_{Y_{\C}}$, so is $z-w_{\C}^{\sigma}$ 
for all $\sigma\in \Aut(\C/\Qb)$.
In particular, we have 
\begin{equation}\label{sigma323}
\rho_{Y_{\C}}(w_{\C})=\rho_{Y_{\C}}(w_{\C}^{\sigma})
\quad\text{ for all $\sigma\in \Aut(\C/\Qb)$.}
\end{equation}
Let $w_S\in \CH^2(Y\times S)_{\hom}$ be a model of the cycle $w_{\C}$.
The cycle $w_S$ 
defines a normal function $\nu_S:S(\C)\ra J^2(Y_{\C})\times S(\C)$,
$s\mapsto (\rho_{Y_{\C}}(w_{S,s}),s)$ where $w_{S,s}\in \CH^2(Y_{\C})_{\hom}$ 
denotes the section of $w_S$ for a closed point $s\in S(\C)$.
Then, \eqref{sigma323} means that the normal function $\nu_S$ is constant on
the subset $S(\C)\setminus S(\Qb)$. 
Since the set $S(\C)\setminus S(\Qb)$ is analytically dense in $S(\C)$,
we have $\nu_S$ is constant on $S(\C)$. 
Let $s\in S(\Qb)$ be any $\Qb$-point, and take
the section $w_{S,s}\in \CH^2(Y)_{\hom}$. 
By the above, $\rho_{Y_{\C}}(w_{S,s})=\rho_{Y_{\C}}(w_{\C})$.  
Clearly $w_{S,s}$ has a support on $D$, so it is the desired cycle.

We thus complete the proof of Theorem~\ref{asakura2}.
\end{pf}
Finally we remark:
\begin{thm}[\cite{masanori} Theorem~4.9]\label{asakura3}
The conjecture~$\ref{mainconj}$ implies
the Beilinson conjecture~$\ref{Beilinson}$~$(1)$.
In particular, it implies the finiteness of the rank of $\CH_0(X)$
for any projective nonsingular variety $X$ over a number field.
\end{thm}
\begin{pf}
Put $\MM_{\Qb}=\MHM(\Spec\Qb)$.
Then the higher Abel-Jacobi map 
$\rho^{\nu}_X$
factors as follows:

\setlength{\unitlength}{1mm}
\begin{picture}(155,37)(-76,-27)
\put(-42,-23){$\Ext^{\nu}_{\MM_{\Qb}}(\Q(0),H^{2r-\nu}(X\ot\Qb)(r))$.}
\put(-36,-3){\vector(1,-1){14}}
\put(-25,2){\vector(1,0){22}}
\put(-5,-17){\vector(1,1){14}}
\put(-15,5){$\rho^{\nu}_X$}
\put(-50,1){$\Gr^{\nu}_F\CH^r(X)$}
\put(0,1){$\Ext^{\nu}_{\MM(\C)}(\Q(0),H^{2r-\nu}(X\ot\C)(r))$}
\end{picture}
By Proposition~\ref{carl325}~(2), the extension groups in $\MM_{\Qb}$
of degree $\geq2$ vanishes.
Therefore we have $F^2\CH^r(X)=0$ by Conjecture \ref{mainconj}.

The latter assertion follows from the Mordell-Weil theorem because of
the isomorphism
$\CH_0(X)_{\dego}\os{\sim}{\ra}\Alb(X)(k)$.
\end{pf}

\medskip
\begin{flushleft}
\noindent
Research Institute for Mathematics Sciences,
Kyoto University,
Oiwakecho, Sakyo-ku,
Kyoto, 606-8502, JAPAN

\medskip

\noindent
E-mail address : asakura@@kurims.kyoto-u.ac.jp
\end{flushleft}
\end{document}